\documentclass[9pt]{amsart}
\usepackage{amsfonts}
\usepackage{graphicx}
\usepackage{amscd}
\usepackage{rotating}
\usepackage{amsmath,amsfonts,amssymb,amsthm,mathrsfs,xypic}
\xyoption{curve}
\usepackage{fancybox}

\newtheorem{thm}{Theorem}[section]

\newtheorem{lem}[thm]{Lemma}

\newtheorem{prop}[thm]{Proposition}
\newtheorem{rem}[thm]{Remark}

\numberwithin{equation}{section}

\newcommand{\s}{\hfill\blacksquare}

\newcommand{\Cok}{\operatorname{Coker}}
\newcommand{\Hom}{\operatorname{Hom}}

\begin{document}
\title [Unbounded ladders induced by Gorenstein algebras]{Unbounded ladders induced by Gorenstein algebras}
\author [P. Zhang, Y. H. Zhang, G. D. Zhou, L. Zhu]{Pu Zhang, Yuehui Zhang, Guodong Zhou, Lin Zhu}
\thanks{\it 2010 Mathematical Subject Classification. \ 18E30, 16E35, 18A40, 18A22, 16G10.}
\thanks{Keywords: periodic ladder, Serre functor, Gorenstein algebra, Gorenstein-projective module, splitting recollement}
\thanks{Supported by the NSFC 11271251 and 11431010.}

\begin{abstract} The derived category $D({\rm Mod}A)$ of a
Gorenstein triangular matrix algebra $A$ admits an unbounded ladder;
and this ladder restricts to $D^-({\rm Mod})$ {\rm(}resp. $D^b({\rm
Mod})$, $D^b({\rm mod})$, $K^b({\rm proj})${\rm)}. A left
recollement of triangulated categories with Serre functors sits in a
ladder of period $1$; as an application, the singularity category of
$A$ admits a ladder of period $1$.
\end{abstract}
\maketitle
Recollements  ([BBD]) provide a
powerful tool for studying problems in triangulated
categories and algebraic geometry. To study mixed categories, ladders have been introduced
([BGS], [AHKLY]). Recollements are ladders of height $1$; while ladders of
height $\ge 2$ give more information ([AHKLY], [HQ]). A fundamental question
is when unbounded ladders occur naturally in representation
theory. This essentially deals with the existence of infinite adjoint sequences.
It is known that if $A$ is an algebra of finite global dimension, then any recollement of derived category $D({\rm Mod} A)$ sits in an unbounded ladder ([AHKLY, 3.7]).

Let $A: = \left(\begin{smallmatrix}B&0\\M&C\end{smallmatrix}\right)$
with $M$ a $C$-$B$-bimodule. An algebra $\Lambda$ is of this form if
and only if $\Lambda$ has an idempotent $e$ with $(1-e)\Lambda e=0$.
It is well-known that there is a recollement $(D({\rm Mod} B), \
D({\rm Mod} A), \ D({\rm Mod} C))$ of derived categories ([AHKL],
[Han]); and in fact, it sits in a ladder of height $2$ ([AHKLY,
3.4]). We further claim that it sits in an unbounded ladder,
provided that $A$, $B$ and $C$ are Gorenstein algebras. This
unbounded ladder enjoys pleasant properties in the sense that it
restricts to $D^-({\rm Mod})$ {\rm(}resp. $D^b({\rm Mod})$,
$D^b({\rm mod})$, $K^b({\rm proj})${\rm)}.

For an adjoint pair $(F, G)$ of categories with Serre functors, $F$ (resp. $G$) always has a left (resp. right) adjoint.
So a left recollement of triangulated categories with Serre functors sits in a ladder of period $1$. As an application, the singularity
category ([B], [O]) of a Gorenstein triangular matrix algebra
admits a ladder of period $1$ (Thm. \ref{matrixalg1}), via the stable category of Gorenstein-projective modules
([EJ], [B], [H2]).
\section {\bf Preliminaries}
\subsection{} Let $\mathcal C',$ $\mathcal C$ and $\mathcal C''$  be triangulated
categories. {\it A recollement} $(\mathcal C', \mathcal C, \mathcal
C'', i^*, i_*, i^!, j_!, j^*, j_*)$ of $\mathcal C$ relative
to $\mathcal C'$ and $\mathcal C''$ ([BBD]) is a diagram of triangle
functors
\begin{center}
\begin{picture}(100,36)
\put(2,22){\makebox(3,1){${\mathcal C}'$}}
\put(40,30){\vector(-1,0){30}}
\put(10,23){\vector(1,0){30}}
\put(40,15){\vector(-1,0){30}}
\put(47,22){\makebox(3,1){${\mathcal
C}$}}
\put(85,30){\vector(-1,0){30}}
\put(55,23){\vector(1,0){30}}
\put(85,15){\vector(-1,0){30}}
\put(93,22){\makebox(3,1){${\mathcal
C}''$}}
\put(25,33){\makebox(3,1){\scriptsize$i^*$}}
\put(25,26){\makebox(3,1){\scriptsize$i_*$}}
\put(25,18){\makebox(3,1){\scriptsize$i^!$}}
\put(71,33){\makebox(3,1){\scriptsize$j_!$}}
\put(71,26){\makebox(3,1){\scriptsize$j^*$}}
\put(71,18){\makebox(3,1){\scriptsize$j_*$}}
\put(220,22){$(1.1)$}
\end{picture}
\end{center}
\vskip-20pt \noindent satisfying the following conditions:

(R1) \ \ $(i^*, i_*), (i_*, i^!), (j_!, j^*)$ and $(j^*, j_*)$ are
adjoint pairs;

(R2) \ \ $i_*, j_!$ and $j_*$ are fully faithful;

(R3) \ \ $j^*i_*=0$ {\rm(}and thus $i^*j_! = 0 = i^!j_*${\rm)};

(R4) \ \ for $X\in\mathcal C$ there are distinguished triangles \
$j_!j^*X \stackrel{\epsilon_X}\longrightarrow X
\stackrel{\eta_X}\longrightarrow i_*i^*X\rightarrow
(j_!j^*X)[1]$ and \ $i_* i^!X\stackrel{\omega_X}\longrightarrow
X\stackrel{\zeta_X}\longrightarrow j_*j^*X\rightarrow (i_*
i^!X)[1],$ where the marked morphisms are the counits and the units of the adjunctions.

\vskip5pt

{\it A left {\rm (resp.} right{\rm)} recollement}  of
$\mathcal C$ relative to $\mathcal C'$ and $\mathcal C''$ is
the upper (resp. lower) two rows of $(1.1)$ satisfying the same
conditions involving only these functors ([P], [K\"o]. For other or related names see e.g. [BGS], [M], [BO], [Kr], [IKM]).  {\it An opposed
recollement} of $\mathcal C$ relative to $\mathcal C'$ and $\mathcal
C''$ is a diagram

\begin{center}
\begin{picture}(100,40)
\put(0,20){\makebox(3,1){${\mathcal C}'$}}
\put(10,27){\vector(1,0){30}}

\put(40,20){\vector(-1,0){30}}

\put(10,13){\vector(1,0){30}}

\put(50,20){\makebox(3,1){${\mathcal C}$}}

\put(60,27){\vector(1,0){30}}

\put(90,20){\vector(-1,0){30}}

\put(60,13){\vector(1,0){30}}

\put(100,20){\makebox(3,1){${\mathcal C}''$}}

\put(25,30){\makebox(3,1){\scriptsize$i_{-1}$}}
\put(25,23){\makebox(3,1){\scriptsize$j_0$}}
\put(25,16){\makebox(3,1){\scriptsize$i_1$}}
\put(75,30){\makebox(3,1){\scriptsize$j_{-1}$}}
\put(75,23){\makebox(3,1){\scriptsize$i_0$}}
\put(75,16){\makebox(3,1){\scriptsize$j_1$}}
\put(275,15){\makebox(25,1)} \put(220,15){}
\end{picture}
\end{center}
\vskip-10pt\noindent such that $(\mathcal C'', \mathcal C, \mathcal C', j_{-1},
i_0, j_1, i_{-1}, j_0, i_1)$ is a recollement of $\mathcal C$
relative to $\mathcal C''$ and $\mathcal C'.$

\begin{lem}\label{criterionupperrec} $(1)$  Given the upper two rows
of  $(1.1)$, the following are equivalent$:$

$({\rm i})$ \  it is a left recollement$;$

$({\rm ii})$ \ $(i^*, i_*)$ and $(j_!, j^*)$ are adjoint pairs,
$i_*$ and $j_!$ are fully faithful, and ${\rm Im}i_{*}={\rm
Ker}j^{*};$

$({\rm iii})$ \   $(i^*, i_*)$ and $(j_!, j^*)$ are adjoint pairs,
$i_*$ and $j_!$ are fully faithful,  and ${\rm Im}j_{!}={\rm
Ker}i^{*}$.
\vskip5pt

$(2)$ \ {\rm(see e.g. [IKM, 1.7])} \ Given diagram $(1.1)$ of
triangle functors, the following are equivalent$:$

$({\rm i})$ \  it is a recollement$;$

$({\rm ii})$ \  it satisfies ${\rm(R1)}$ and ${\rm(R2)}$,  and ${\rm
Im}i_{*}={\rm Ker}j^{*}, \ {\rm
Im}j_{!}={\rm Ker}i^{*}$ and ${\rm
Im}j_* ={\rm Ker}i^!;$

$({\rm iii})$ \  it satisfies ${\rm(R1)}$ and ${\rm(R2)}$, and any one of the equalities in $(2)$$({\rm ii})$.
\end{lem}

\subsection{} \ {\it A ladder} ({\rm [BGS, 1.2], [AHKLY, Sect. 3]})
 is a finite or an infinite
diagram of triangle functors:
\begin{center}
\begin{picture}(120,60)

\put(0,20){\makebox(3,1){${\mathcal C}'$}}

\put(25,45){\vdots}

\put(10,35){\vector(1,0){30}}

\put(40,27){\vector(-1,0){30}}

\put(10,20){\vector(1,0){30}}

\put(40,13){\vector(-1,0){30}}

\put(10,5){\vector(1,0){30}}

\put(25,-10){\vdots}

\put(50,20){\makebox(3,1){${\mathcal C}$}}

\put(75,45){\vdots}

\put(60,35){\vector(1,0){30}}

\put(90,27){\vector(-1,0){30}}

\put(60,20){\vector(1,0){30}}

\put(90,13){\vector(-1,0){30}}
\put(60,5){\vector(1,0){30}}
\put(75,-10){\vdots}

\put(100,20){\makebox(3,1){${\mathcal
C}''$}}

\put(23,39){\makebox(3,1){\scriptsize$i_{-2}$}}
\put(75,39){\makebox(3,1){\scriptsize$j_{-2}$}}
\put(75,30){\makebox(3,1){\scriptsize$i_{-1}$}}
\put(23,30){\makebox(3,1){\scriptsize$j_{-1}$}}
\put(23,23){\makebox(3,1){\scriptsize$i_0$}}
\put(75,23){\makebox(3,1){\scriptsize$j_0$}}
\put(23,16){\makebox(3,1){\scriptsize$j_1$}}
\put(75,16){\makebox(3,1){\scriptsize$i_1$}}
\put(23,8){\makebox(3,1){\scriptsize$i_2$}}
\put(75,8){\makebox(3,1){\scriptsize$j_2$}}
\put(275,15){\makebox(25,1)} \put(220,15){\rm (1.2)}
\end{picture}
\end{center}

\vskip10pt \noindent such that any two consecutive rows form a left
or right recollement {\rm(}or equivalently, any three consecutive
rows form a recollement or an opposed recollement) of $\mathcal C$
relative to $\mathcal C'$ and $\mathcal C''$. Its {\it height} is
the number of rows minus $2$. Ladders of height $0$ (resp. $1$) are exactly left
or right recollements (resp. recollements or opposed recollements). A ladder is {\it unbounded}
 if it goes infinitely both upwards and downwards.

\vskip5pt

A two-sided infinite sequence \ $(\cdots, F_{-1}, F_0, F_1, \cdots)$ \ of additive functors is {\it an infinite adjoint sequence}, if
$(F_n, F_{n+1})$ is an adjoint pair for each $n\in\Bbb Z$. In such a sequence if some $F_i$ is a triangle functor then so are all $F_n$'s (see e.g. [Ke1, 6.7]).

\begin{lem}\label{laddercriterion}  Recollement ${\rm (1.1)}$ sits in an unbounded ladder if and only if there is an infinite adjoint sequence \ ${\rm(}\cdots, \ F_{-1}, \ i^*, \ i_*, \ i^!, \ F_1, \ \cdots{\rm)}$.
\end{lem}
\subsection{} {\it An equivalence} of left recollements ([PS, 2.5], [FP]) is a triple $(F', F, F'')$ of triangle-equivalences such
that

\begin{center}
\begin{picture}(100,65)
\put(0,60){\makebox(3,1){${\mathcal C}'$}}
\put(35,65){\vector(-1,0){25}}
\put(10,57){\vector(1,0){25}}
\put(40,60){\makebox(3,1){${\mathcal C}$}}
\put(75,65){\vector(-1,0){25}} \put(50,57){\vector(1,0){25}}

\put(85,60){\makebox(3,1){${\mathcal C}^{\prime\prime}$}}

\put(21,68){\makebox(3,1){\scriptsize$i^\ast$}}
\put(21,60){\makebox(3,1){\scriptsize$i_\ast$}}

\put(63,68){\makebox(3,1){\scriptsize$j_!$}}
\put(63,60){\makebox(3,1){\scriptsize$j^\ast$}}

\put(-12,38){$F'$}

\put(45,38){$F$}

\put(90,38){$F''$}

\put(0,54){\vector(0,-1){26}}
\put(40,54){\vector(0,-1){26}}
\put(85,54){\vector(0,-1){26}}

\put(0,21){\makebox(3,1){${\mathcal D}'$}}

\put(35,26){\vector(-1,0){25}}
\put(10,16){\vector(1,0){25}}

\put(40,21){\makebox(3,1){${\mathcal D}$}}

\put(75,26){\vector(-1,0){25}}
\put(50,16){\vector(1,0){25}}
\put(85,21){\makebox(3,1){${\mathcal D}^{\prime\prime}$}}

\put(21,31){\makebox(3,1){\scriptsize$i_\mathcal D^\ast$}}
\put(21,20){\makebox(3,1){\scriptsize$i^\mathcal D_\ast$}}
\put(63,31){\makebox(3,1){\scriptsize$j^\mathcal D_!$}}
\put(63,20){\makebox(3,1){\scriptsize$j_\mathcal D^\ast$}}
\end{picture}
\end{center}
\vskip-15ptcommutes. Similarly we have an equivalence of (right,
opposed) recollements.
\vskip5pt
We call  $(\mathcal C', \mathcal C, \mathcal C'',
j_{2t-1}, i_{2t}, j_{2t+1}, i_{2t-1}, j_{2t}, i_{2t+1})$ in ladder
$(1.2)$ {\it the $t$-th recollement}, \ $(\mathcal C', \mathcal C,
\mathcal C''$, $i_{2t}, j_{2t+1},$ $ i_{2t+2}, j_{2t}, i_{2t+1},
j_{2t+2})$  {\it the $t$-th opposed recollement}, and the left
(right) recollement sitting in the $t$-th recollement {\it the $t$-th
left {\rm(}right{\rm)} recollement}. An unbounded ladder  $(1.2)$ is {\it periodic}, if
there is an integer $t\ge 1$ such that the $t$-th left
recollement is equivalent to the $0$-th one. Such a
minimal $t$ is called {\it the period}. The following describes the period via the associated TTF tuple, and justifies the terminology.
\begin {lem} \label{periodicladder} \  $(1)$  \ Given recollements
$(\mathcal C', \mathcal C, \mathcal C'')$ and $(\mathcal D',
\mathcal D, \mathcal D'')$, the following  are equivalent$:$

${\rm (i)}$ \ \ they are equivalent$;$

${\rm (ii)}$ \  there is a triangle-equivalence $F: \mathcal
C\rightarrow \mathcal D$ such that $F({\rm Im}j_!) = {\rm
Im}j_!^\mathcal D$,  $F({\rm Im}i_*) = {\rm Im}i_*^\mathcal D$ and
$F({\rm Im}j_*) = {\rm Im}j_*^\mathcal D;$

${\rm (iii)}$ \ \ there is a triangle-equivalence $F: \mathcal
C\rightarrow \mathcal D$ such that one of the equalities in ${\rm (ii)}$ holds.

\vskip5pt

$(2)$ \ Given a ladder of period $t$, then the $(qt+l)$-th {\rm(}left, right, opposed{\rm)}
recollement  is equivalent to the $l$-th {\rm(}left, right,
opposed{\rm)} recollement for $q\in\Bbb Z$ and $l = 0, \cdots, t-1$,
under the same equivalence.
\vskip5pt

$(3)$ \ Given an unbounded ladder ${\rm(1.2)}$, the following are equivalent$:$

${\rm (i)}$ \ \ it is of period $t;$

${\rm (ii)}$ \ $t$ is the minimal positive integer such that there is a
triangle-equivalence $F: \mathcal C\rightarrow \mathcal C$ satisfying
$F({\rm Im}i_{2t+1}) = {\rm Im}i_1, \ \ F({\rm Im}i_{2t}) = {\rm
Im}i_0$ and $F({\rm Im}i_{2t-1}) = {\rm Im}i_{-1};$

${\rm (iii)}$ \  $t$ is the minimal positive integer such that there is a
triangle-equivalence $F: \mathcal C\rightarrow \mathcal C$ satisfying
one of the equalities in ${\rm (ii)}$.
\end{lem}
\subsection{} If no specified, modules are right modules.
For algebra $A$ over a field, denote by Mod$A$ (resp. $A$-Mod) the
category of right (resp. left) $A$-modules. If $A$ is
finite-dimensional, then we denote by mod$A$ (resp. $A$-mod) the
category of finitely generated right (resp. left) $A$-modules, and
by $\mathcal{GP}(A)$ the full subcategory of mod$A$ consisting of
Gorenstein-projevtive modules ([EJ]). Then $\mathcal{GP}(A)$ is a
Frobenius category whose projective-injective objects are exactly
projective modules ([Be]), and hence the stable category
$\underline{\mathcal {GP}(\mathcal A)}$ is triangulated ([H1]). A
finite-dimensional algebra $A$ is {\it Gorenstein} if ${\rm inj.dim}
_AA< \infty$ and ${\rm inj.dim} A_A < \infty$.

\vskip5pt

Let $K^{b}({\rm proj} A)$ (resp. $K^{b}({\rm inj} A)$) be the
homotopy category of bounded complexes of finitely generated
projective (resp. injective) right $A$-modules,  $D({\rm Mod}A)$
(resp. $D^-({\rm Mod}A)$, $D^b({\rm Mod}A)$) the unbounded (resp.
upper bounded, bounded) derived category of Mod$A$, and $D^{b}({\rm
mod}A)$ the bounded derived category of mod$A$. Note that $D({\rm
Mod}A)$ is compactly generated by $A_A$ (see [S]; also [BN]).

\vskip5pt

For a triangulated category $\mathcal T$ with coproducts, denote by
$\mathcal T^c$ the full subcategory of $\mathcal T$ consisting of
compact objects. Then $D^c({\rm Mod}A) = K^{b}({\rm proj} A)$
([N1]).

\section {\bf Main results}
\begin{thm} \label{matrixalg2} \ Let $B$ and $C$ be
Gorenstein algebras and $_CM_B$ a $C$-$B$-bimodule, such that $A =
\left(\begin{smallmatrix}B & 0 \\ M & C\end{smallmatrix}\right)$ is
Gorenstein. Then there is an unbounded ladder $(D({\rm Mod} B), \
D({\rm Mod} A), \ D({\rm Mod} C))$ of derived categories.
\end{thm}
\begin{rem} For the Gorensteinness of $A:=\left(\begin{smallmatrix}B & 0 \\ M & C\end{smallmatrix}\right)$ we refer to {\rm[C]} and {\rm[Z, Thm. 2.2]}. If $B$ and $C$ are Gorenstein, then $A$ is Gorenstein if and only if ${\rm
proj.dim}_CM$ and ${\rm proj.dim}M_B$ are finite {\rm([C,
Thm. 3.3])}. Also note that ${\rm gl.dim}A \ge {\rm max}\{{\rm gl.dim}B, \ {\rm gl.dim}C\}$.

For example, let $A$ be the algebra given by quiver
\xymatrix{{\bullet}\ar@(ur,ul)^{\lambda_3}\ar[r]^\beta &
\ar@(ur,ul)^{\lambda_1}\bullet
&\ar[l]_\alpha\ar@(ur,ul)^{\lambda_2}\bullet} \ \ and relations  \
$\lambda_1^2, \ \lambda_2^2, \ \lambda_3^2, \ \alpha\lambda_2
-\lambda_1\alpha, \ \beta\lambda_3 -\lambda_1\beta$. Then $A
=\left(\begin{smallmatrix} B & 0 \\ _CM_B &
C\end{smallmatrix}\right)
= \left(\begin{smallmatrix} C & 0& 0 \\ C & C & 0 \\
0& C & C\end{smallmatrix}\right)$,  where $C: = k[x]/\langle
x^2\rangle$, $B: = T_2(C): = \left(\begin{smallmatrix}C & 0 \\ C &
C\end{smallmatrix}\right)$ and $_CM_B: = \ _C(0, C)_{T_2(C)}$. Since
${\rm proj.dim}_CM = 0$ and ${\rm proj.dim}M_{T_2(C)}=1$,  $A$ is
Gorenstein of ${\rm gl.dim} A = \infty.$\end{rem}

\subsection{} Let $A =
\left(\begin{smallmatrix}B & 0 \\ _CM_B &
C\end{smallmatrix}\right)$. A right $A$-module is given by $(X_B,
Y_C)_\phi,$ where $X_B\in$ mod $B$, $Y_C\in$ mod $C$, and $\phi: \
Y\otimes_CM_B\rightarrow X_B$ is a right $B$-map. A right $A$-map
$(X_B, Y_C)_\phi\rightarrow (X'_B, Y'_C)_{\phi'}$ is given by $(f,
g)$ with $f\in {\Hom}_B (X_B, X'_B)$ and $g\in \Hom_C(Y_C, Y'_C)$,
such that \ $f\phi = \phi'(g\otimes_B {\rm Id}_M)$. A left
$A$-module is given by $\left(\begin{smallmatrix}
_BX\\
_CY
\end{smallmatrix}\right)_\phi$,  where $_BX\in B$-mod, \ $_CY\in C$-mod,
and $\phi: \ _CM\otimes_B X\rightarrow \ _CY$ is a left $C$-map. A
left $A$-map $\left(\begin{smallmatrix}
_BX\\
_CY
\end{smallmatrix}\right)_\phi\rightarrow \left(\begin{smallmatrix}
_BX'\\
_CY'
\end{smallmatrix}\right)_{\phi'}$
is given by $(f, g)$ with $f\in {\Hom}_B (_BX, \ _BX')$ and $g\in
\Hom_C(_CY, \ _CY')$, such that \ $g\phi = \phi'({\rm Id}_M\otimes_B
f)$. The projective right $A$-modules are exactly \ $(P_B, 0)$ and
$(Q\otimes_CM, Q_C)_{\rm Id}$,  where $P_B\in {\rm proj}B$ and
$Q_C\in {\rm proj}C.$ The projective left $A$-modules are exactly
$\left(\begin{smallmatrix}
   _BP \\
   M\otimes_BP
\end{smallmatrix}\right)_{\rm Id}
$ and $\left( \begin{smallmatrix}
   0\\
   _CQ
\end{smallmatrix}\right)$, where $_BP\in B\mbox{-}{\rm proj}$ and $_CQ\in C\mbox{-}{\rm proj}.$
See [ARS, p.73].

\subsection{} Let $A$ be an algebra over a field with idempotent $e$. The ideal $AeA$ is {\it stratifying} ([CPS, 2.1.1]),
if the multiplication map $m: Ae\otimes_{eAe}eA \rightarrow AeA$ is
injective and ${\rm Tor}^n_{eAe}(Ae, eA) = 0$ for $n\ge 1$. As
pointed out by S. K\"onig and H. Nagase [KN, Rem. 3.2], $_A(AeA)$ \
(resp. $(AeA)_A$) \ is projective if and only if $_{eAe}(eA)$ \
(resp. $(Ae)_{eAe}$) \ is projective and the map $m$ is injective.
Thus, if $AeA$ is projective either as a left or a right $A$-module,
then $AeA$ is a stratifying ideal.

\begin{lem} \label{cps} \ {\rm (see e.g. [AHKL, 4.5], [Han, 2.1])} \ If $AeA$ is
a stratifying ideal, then there is a recollement
\begin{center}
\begin{picture}(120,38)
\put(-15,20){\makebox(-22,1){$D({\rm Mod}A/AeA)$}}
\put(40,28){\vector(-1,0){30}} \put(10,20){\vector(1,0){30}}
\put(40,13){\vector(-1,0){30}} \put(50,20){\makebox(25,0.8){$D({\rm
Mod}A)$}} \put(115,28){\vector(-1,0){30}}
\put(85,20){\vector(1,0){30}} \put(115,13){\vector(-1,0){30}}
\put(132,20){\makebox(27,0.5){$D({\rm Mod} \ eAe)$}}
\put(25,31){\makebox(3,1){\scriptsize$i^\ast$}}
\put(25,23){\makebox(3,1){\scriptsize$i_\ast$}}
\put(25,16){\makebox(3,1){\scriptsize$i^!$}}
\put(100,31){\makebox(3,1){\scriptsize$j_!$}}
\put(100,23){\makebox(3,1){\scriptsize$j^\ast$}}
\put(100,16){\makebox(3,1){\scriptsize$j_\ast$}}
\end{picture}
\end{center}
\vskip-10pt \noindent where\vskip-10pt \noindent
\begin{align*} & i^* = -\overset{\rm L}\otimes_A A/AeA,  \ \ \ \ \ \ \ \  i_* = -\overset{\rm L}\otimes_{A/AeA}A/AeA, \ \ \ \ \ \ \ \  i^! = {\rm
R}\Hom_A(A/AeA, -), \\& j_! = -\overset{\rm L}\otimes_{eAe} eA, \   \ \ \ \ \ \ \
\ \ \  j^* = -\overset{\rm L}\otimes_{A} Ae,  \ \ \ \ \ \ \ \ \ \ \ \ \ \ \ \ \ \
\ j_* = {\rm R}\Hom_{eAe}(Ae, -).
\end{align*}
\end{lem}

\subsection{} Let $\mathcal T$ be a triangulated
category compactly generated by $\mathcal S_0$. Denote by $\langle
\mathcal S_0\rangle$ the smallest triangulated subcategory of
$\mathcal T$ containing $\mathcal S_0$ and closed under coproducts.
Brown representability  (A. Neeman [N2, Thm. 3.1]) claims that every
cohomological functor $F: \mathcal T^{op} \rightarrow {\rm Ab}$
which sends coproducts to products is representable {\rm(}i.e.,
$F\cong \Hom_\mathcal T(-, X)$ for some $X\in\mathcal T${\rm)}, and
that \ $\mathcal T = \langle \mathcal S_0\rangle$.  And, Brown
representability for the dual (H. Krause [Kr, Thm. A]) claims that \
$\mathcal T$ has products, and that every cohomological functor $F:
\mathcal T \rightarrow {\rm Ab}$ which sends products to products is
representable {\rm(}i.e., $F\cong \Hom_\mathcal T(X, -)$ for some
$X\in\mathcal T${\rm)}.

\vskip5pt

Using Brown representability one has

\begin{lem}\label{GhasrightadjandFkeepcompact} \ {\rm ([N2, Thm. 4.1 and 5.1])} \ Let $F: \mathcal C
\rightarrow \mathcal D$ be a triangle functor between compactly
generated triangulated categories, with a right adjoint $G$. Then
the following are equivalent$:$

${\rm(i)}$ \ $G$ admits a right adjoint$;$

${\rm(ii)}$ \ $F$ preserves compact objects$;$

${\rm(iii)}$ \ $G$ preserves coproducts.
\end{lem}

Using Brown representability for the dual one has

\begin{lem}\label{leftadjandproduct} \ Let $F: \mathcal C
\rightarrow \mathcal D$ be a triangle functor between triangulated
categories, where $\mathcal C$ is compactly generated. Then $F$
admits a left adjoint if and only if $F$ preserves products {\rm(}we
are not assuming that $\mathcal D$ has products{\rm)}.
\end{lem}
\noindent{\bf Proof.} \ The ``only if" part is well-known. For the
``if" part, applying Brown representability for the dual to functor
$\Hom_\mathcal D(Y, F-): \mathcal C \rightarrow {\rm Ab}$, for each
object $Y\in\mathcal D$, we then see that $F$ admits a left adjoint.
$\s$

\vskip5pt

We need the following result due to P. Balmer, I. Dell'ambrogio and
B. Sanders.

\begin{lem}\label{preservesproduct} \ {\rm ([BDS, Lemma
2.6(b)])} \ Let $F: \mathcal C \rightarrow \mathcal D$ be a triangle
functor between compactly generated triangulated categories, with a
right adjoint $G$. Assume that $F$ preserves compacts, and the
restriction $F|_{\mathcal C^c}: \mathcal C^c \rightarrow \mathcal
D^c$ admits a left adjoint. Then $F$ preserves products.
\end{lem}

\subsection{} Let $A$ be a finite-dimensional algebra over field $k$. Using a
hoprojective (resp. hoinjective) resolution of a complex in $D({\rm
Mod}A)$ ([S]; [BN]) one has the  characterizations:
$$K^b({\rm proj}A) = \{P\in D({\rm Mod}A) \ | \ {\rm dim}_k(\bigoplus\limits_{i\in\Bbb Z} \Hom_{D({\rm Mod} A)}(P, Y[i])) < \infty, \ \forall \ Y\in D^b({\rm mod}A) \},$$
and
$$K^b({\rm inj}A) = \{I\in D({\rm Mod}A) \ | \ {\rm dim}_k(\bigoplus\limits_{i\in\Bbb Z} \Hom_{D({\rm Mod} A)}(Y[i], I)) < \infty, \ \forall \ Y\in D^b({\rm mod}A) \}.$$
One has also the  characterization:
$$D^b({\rm mod}A) = \{X\in D({\rm Mod}A) \ | \ {\rm dim}_k(\bigoplus\limits_{i\in\Bbb Z} \Hom_{D({\rm Mod} A)}(P, X[i])) < \infty, \ \forall \ P\in K^b({\rm proj}A) \}.$$
See L. Angeleri H\"ugel, S. K\"onig, Q. H. Liu and D. Yang [AHKLY,
Lemma 2.4].  Using these one has

\begin{lem}\label{FkeepcompactGkeepDbmod}  \ Let $A$ and $B$ be
finite-dimensional algebras, and $F: D({\rm Mod} A) \rightarrow
D({\rm Mod} B)$ a triangle functor with a right adjoint $G$. Then

${\rm(i)}$ \ {\rm ([AHKLY, Lemma 2.7])} \ $F$ preserves $K^b({\rm
proj})$ if and only if $G$ preserves $D^b({\rm mod})$.

${\rm(ii)}$ \ {\rm ([HQ, Lemma 1])} \ $F$ preserves $D^b({\rm mod})$
if and only if $G$ preserves $K^b({\rm inj})$.
\end{lem}

\subsection {} Let $\mathcal C$ be a Hom-finite category over field $k$. A $k$-linear functor $S:
\mathcal C \rightarrow \mathcal C$ is {\it a right Serre functor},
if for any objects $X$ and $Y$ there is a $k$-isomorphism
$\Hom_\mathcal C(X, Y)\cong \Hom_\mathcal C(Y, SX)^*$ which is
natural in $X$ and $Y$, where $(-)^* = \Hom_k(-, k)$. We say that
$\mathcal C$ {\it has a Serre functor} if $\mathcal C$ has a right
Serre functor which is an equivalence, or equivalently, $\mathcal C$
has both a right and left Serre functor ([BK], [RV]). If $\mathcal
C$ is a Hom-finite Krull-Schmidt triangulated category over an
algebraically closed field $k$, then $\mathcal C$ has a Serre
functor if and only if $\mathcal C$ has Auslander-Reiten triangles
(note that the assumption that $k$ is algebraically closed is only
used in the ``only if" part. See I. Reiten and M. Van den Bergh [RV,
Thm. 2.4]).

\vskip5pt

The following observation will play an important role in this paper.

\begin {lem} \label{withserre} \ Let $\mathcal C$ and $\mathcal D$
be categories with Serre functors, $F: \mathcal{C}\rightarrow
\mathcal {D}$ an additive functors with a right adjoint $G$. Then
$F$ admits a left adjoint $S_\mathcal {C}^{-1}GS_\mathcal {D},$  and
$G$ admits a right adjoint $S_\mathcal {D}FS_\mathcal {C}^{-1}$,
where $S_\mathcal C$ and $S_\mathcal D$ are the right Serre functors
of $\mathcal C$ and $\mathcal D$, respectively. \end{lem}
\noindent{\bf Proof.} \ For $X\in \mathcal C$ and $Y\in\mathcal D$
we have
$$\Hom_\mathcal C(S_\mathcal {C}^{-1}GS_\mathcal {D}Y, X) \cong
\Hom_\mathcal C(X, GS_\mathcal {D}Y)^* \cong \Hom_\mathcal D(FX,
S_\mathcal D Y)^* \cong \Hom_\mathcal D(Y, FX).$$ Similarly $(G,
S_\mathcal {D}FS_\mathcal {C}^{-1})$ is an adjoint pair. $\s$

\vskip5pt

We also need the following result due to D. Happel.

\begin{lem}\label{kbp=kbj} {\rm ([Hap2, Lemma 1.5, Thm.
3.4])} \ \ Let $A$ be a finite-dimensional algebra. Then $A$ is
Gorenstein if and only if  $K^b({\rm proj}A) = K^b({\rm inj}A)$ in
$D^b({\rm mod} A)$. In this case $K^b({\rm proj}A)$ has a Serre
functor.
\end{lem}

\subsection{\bf Proof of Theorem \ref{matrixalg2}} \ Put $e: = \left(\begin{smallmatrix}0&0\\ 0&1 \end{smallmatrix}\right)\in A$. Then $AeA = \left(\begin{smallmatrix}0&0\\ M&C \end{smallmatrix}\right)\cong (M, C)$
is a projective right $A$-module, and hence $AeA$ is stratifying. Since
$A/AeA \cong B$ and  \ $eAe\cong C$ as algebras, and
\vskip-10pt$$_A(A/AeA)_B \cong \ _A\left( \begin{smallmatrix}B\\ 0 \end{smallmatrix}\right)_B, \ \ \ _B(A/AeA)_A \cong \ _B(B, 0)_A, \ \ \ _C(eA)_A \cong \ _C(M, C)_A,
\ \ \ _A(Ae)_C \cong \ _A\left( \begin{smallmatrix}0\\ C \end{smallmatrix}\right)_C$$
as bimodules, it follows from Lemma \ref{cps} that there is a recollement
\begin{center}
\begin{picture}(120,40)
\put(-5,20){\makebox(-22,1){$D({\rm Mod}B)$}}
\put(40,28){\vector(-1,0){30}} \put(10,20){\vector(1,0){30}}
\put(40,12){\vector(-1,0){30}} \put(50,20){\makebox(25,0.8){$D({\rm
Mod}A)$}}
\put(115,28){\vector(-1,0){30}}
\put(85,20){\vector(1,0){30}}
\put(115,12){\vector(-1,0){30}}
\put(127,20){\makebox(27,0.5){$D({\rm Mod}C)$}}
\put(25,31){\makebox(3,1){\scriptsize$i^\ast$}}
\put(25,23){\makebox(3,1){\scriptsize$i_\ast$}}
\put(25,15){\makebox(3,1){\scriptsize$i^!$}}
\put(100,31){\makebox(3,1){\scriptsize$j_!$}}
\put(100,23){\makebox(3,1){\scriptsize$j^\ast$}}
\put(100,15){\makebox(3,1){\scriptsize$j_\ast$}}
\put(220,20){$(2.1)$}\end{picture}
\end{center}
\vskip-10pt \noindent where \ \
$i^* = -\overset{\rm L}\otimes_A \left( \begin{smallmatrix}B\\ 0 \end{smallmatrix}\right), \ i_* = -\overset{\rm L}\otimes_B (B,0), \ i^! = {\rm
R}\Hom_A((B,0)_A, -), \ j_! = -\overset{\rm L}\otimes_C(M, C), \  j^* = -\overset{\rm L}\otimes_A\left( \begin{smallmatrix}0\\ C \end{smallmatrix}\right), \ j_* = {\rm R}\Hom_C(\left( \begin{smallmatrix}0\\ C \end{smallmatrix}\right)_C, -).$

\vskip5pt

{\bf Claim 1.} \ There is an infinite sequence $(\cdots, F_{-3}, F_{-2}, F_{-1}, i^*)$ such that any two consecutive functors form an adjoint pair.

\vskip5pt

Since the right adjoint $i_*$ of $i^*$ admits a right adjoint $i^!$,
it follows from Lemma \ref{GhasrightadjandFkeepcompact} that $i^*$
preserves compacts (this could be also seen directly: since $\left(
\begin{smallmatrix}B\\ 0 \end{smallmatrix}\right)$ is projective as a right $B$-module, it follows that $i^* = -\overset{\rm L}\otimes_A \left( \begin{smallmatrix}B\\ 0 \end{smallmatrix}\right)$ preserves
compacts). Since $(B, 0)$ is projective as a right $A$-module, it
follows that $i_* = -\overset{\rm L}\otimes_B (B,0)$ preserves
compacts. Thus $(i^*|_{K^b({\rm proj}A)}, \ i_*|_{K^b({\rm
proj}B)})$ is an adjoint pair. Since $A$ and $B$ are Gorenstein
algebras, by Lemma \ref{kbp=kbj}, $K^b({\rm proj}A)$ and $K^b({\rm
proj}B)$ have Serre functors, and hence $i^*|_{K^b({\rm proj}A)}$
has a left adjoint, by Lemma \ref{withserre}. Applying Lemma
\ref{preservesproduct} to the adjoint pair $(i^*, i_*)$ we know that
$i^*$ preserves products, and hence by Lemma
\ref{leftadjandproduct}, $i^*$ admits a left adjoint, denoted by
$F_{-1}$.

Repeating the above arguments we get {\bf Claim 1}.

\vskip5pt

{\bf Claim 2.} \ There is an infinite sequence $(i^!, G_1, G_2, G_3, \cdots)$ such that any two consecutive functors form an adjoint pair.

\vskip5pt

Since $i_*$ preserves compacts, it follows from Lemma
\ref{GhasrightadjandFkeepcompact} that $i^!$ admits a right adjoint,
denoted by $G_1$.

Since $i^*$ preserves compacts, i.e., $i^*$ preserves $K^b({\rm
proj})$, it follows from Lemma
\ref{FkeepcompactGkeepDbmod}${\rm(i)}$ that $i_*$ preserves
$D^b({\rm mod})$, and hence $i^!$ preserves $K^b({\rm inj})$ by
Lemma \ref{FkeepcompactGkeepDbmod}${\rm(ii)}$. Since we are dealing
with Gorenstein algebras, by Lemma \ref{kbp=kbj} this is exactly to
say that $i^!$ preserves $K^b({\rm proj})$, i.e., $i^!$ preserves
compacts.  It follows from Lemma \ref{GhasrightadjandFkeepcompact}
that $G_1$ admits a right adjoint, denoted by $G_2$.

By the same argument we know that $G_1$ preserves compacts, and
hence by Lemma \ref{GhasrightadjandFkeepcompact},  $G_2$ admits a
right adjoint, denoted by $G_3$. Also, $G_2$ and $G_3$ preserve
compacts. Repeating these arguments we get {\bf Claim 2}.

\vskip5pt

Now Theorem \ref{matrixalg2} follows from Lemma
\ref{laddercriterion}. $\s$

\begin{rem} The unbounded ladder in Theorem \ref{matrixalg2} restricts to $D^-({\rm Mod})$,
$D^b({\rm Mod})$, $D^b({\rm mod})$ and $K^b({\rm proj})$. \ In fact,
since $A$, $B$ and $C$ are Gorenstein, all the functors in
recollement ${\rm(2.1)}$ restrict to $D^-({\rm Mod})$ {\rm(}resp.
$D^b({\rm Mod})$, $D^b({\rm mod})$, $K^b({\rm proj})${\rm)}$;$ then
by Lemmas {\rm\ref{GhasrightadjandFkeepcompact}} and
{\rm\ref{FkeepcompactGkeepDbmod}} we see that all the functors in
the ladder restrict to $K^b({\rm proj})${\rm)} and $D^b({\rm mod})$,
respectively.  By {\rm[AHKLY, Prop. 4.11]} and {\rm[AHKLY, Coroll.
4.9]}, we also see that all the functors in the ladder restrict to
$D^-({\rm Mod})$ and $D^b({\rm Mod})$, respectively.\end{rem}

\section {\bf Ladders of period $1$}
\subsection {} We have
\begin {prop} \label{upperrecollwithserre} \ $(1)$ \ Let $\mathcal C'$, $\mathcal C$ and $\mathcal C''$ be triangulated categories
with Serre functors. Then

${(\rm i)}$ \ Any left {\rm(}right{\rm)} recollement $(\mathcal C',
\mathcal C, \mathcal C'')$ sits in a ladder of period
$1$.

{(\rm ii)} \ Any recollement $(\mathcal C', \mathcal C, \mathcal C'')$
sits in a ladder of period $1$.

\vskip5pt

$(2)$ \ Any
recollement of triangulated category $\mathcal C$ with Serre functor
sits in a ladder of period $1$.
\end{prop} \noindent{\bf
Proof.} \ $(1){(\rm i)}$ \  Let $S_{\mathcal C'}$, $S_{\mathcal C}$ and
$S_{\mathcal C''}$ be the right Serre functors of $\mathcal C',
\mathcal C$ and $\mathcal C''$, respectively. Let
\begin{center}
\begin{picture}(120,36)
\put(12,20){\makebox(-22,1) {$\mathcal C'$}}
\put(40,26){\vector(-1,0){30}}
\put(10,17){\vector(1,0){30}}
\put(38,20){\makebox(25,0.8) {$\mathcal C$}}
\put(90,26){\vector(-1,0){30}}
\put(60,17){\vector(1,0){30}}
\put(90,20){\makebox(25,0.5){$\mathcal C''$}}
\put(25,29){\makebox(3,1){\scriptsize$j_{-1}$}}
\put(25,21){\makebox(3,1){\scriptsize$i_0$}}
\put(75,29){\makebox(3,1){\scriptsize$i_{-1}$}}
\put(75,21){\makebox(3,1){\scriptsize$j_0$}}
\end{picture}
\end{center} \vskip-10pt \noindent be a left recollement. Applying  Lemma
\ref{withserre} to adjoint pair $(j_{-1}, i_0)$  we know that
$j_{-1}$ admits a left adjoint $i_{-2} = S_\mathcal
{C}^{-1}i_0S_\mathcal {C'}: \mathcal C'\rightarrow \mathcal C$, and
that $i_0$ admits a right adjoint $j_1 = S_\mathcal
{C'}j_{-1}S_\mathcal {C}^{-1}: \mathcal C\rightarrow \mathcal C'$.
Similarly,  $i_{-1}$ admits a left adjoint $j_{-2} = S_\mathcal
{C''}^{-1}j_0S_\mathcal {C}$, and $j_0$ admits a right adjoint $i_1
= S_\mathcal {C}i_{-1}S_\mathcal {C''}^{-1}$. By induction we have
\begin{align*} &i_{2n-1} = S_\mathcal {C}^ni_{-1}S^{-n}_\mathcal {C''}: \
\mathcal {C''}\longrightarrow \mathcal {C}, \ \ \ \ \ \  i_{2n} =
S_\mathcal {C}^ni_0S^{-n}_\mathcal {C'}: \ \ \mathcal
{C'}\longrightarrow \mathcal C, \\&
 j_{2n-1} = S_\mathcal {C'}^nj_{-1}S^{-n}_\mathcal {C}: \
\mathcal {C}\longrightarrow \mathcal {C'}, \ \ \ \ \  j_{2n} =
S_\mathcal {C''}^nj_0S^{-n}_\mathcal {C}: \ \mathcal
{C}\longrightarrow \mathcal {C''}.\end{align*} By Lemma
\ref{criterionupperrec}$(2)$ $(\mathcal C', \mathcal C, \mathcal
C'', j_{-1}, i_0, i_1, i_{-1}, j_0, j_1)$ is a recollement, and
hence by Lemma \ref{laddercriterion} we get the desired unbounded
ladder. Since $(S_{\mathcal C'}, S_\mathcal C, S_{\mathcal C''})$ is
an equivalence from the $1$st left recollement to the $0$-th left
recollement, this ladder is of period $1$.

\vskip5pt

${(\rm ii)}$ follows from ${(\rm i)}$ and the fact that one functor in an adjoint
pair uniquely determines another.

\vskip5pt

$(2)$ \ Let $(\mathcal C', \mathcal C, \mathcal C'', i^*, i_*, i^!, j_!, j^*, j_*)$ be a
recollement, and $S$ a right Serre functor of $\mathcal {C}$. Then
$\mathcal C'$ has a right Serre functor $S_\mathcal {C^\prime} =
i^!Si_*$ with $S^{-1}_\mathcal {C^\prime} = i^*S^{-1}i_*;$  and
$\mathcal {C}''$ has a right Serre functor $S_\mathcal
{C^{\prime\prime}} = j^*Sj_!$ with $S^{-1}_\mathcal
{C^{\prime\prime}} = j^*S^{-1}j_*$ (see {\rm P. J$\o$rgensen [J]}. We stress that this
result does not hold for left recollements). Then from $(1){(\rm ii)}$ the assertion follows. $\s$

\subsection{} If $A$ is Gorenstein, then
$\underline {\mathcal {GP}(A)}$ is triangle-equivalent to the
singularity category $D^b({\rm mod}A)/K^b({\rm proj} A)$ ([B,
4.4.1]). So the following gives a ladder of
singularity categories of period $1$.

\vskip5pt

\begin{thm} \label{matrixalg1} \ Let $B$ and $C$ be
Gorenstein algebras and $_CM_B$ a $C$-$B$-bimodule, such that $A =
\left(\begin{smallmatrix}B & 0 \\ M & C\end{smallmatrix}\right)$ is
Gorenstein. Then we have a ladder $(\underline {\mathcal
{GP}(B)}, \ \underline {\mathcal {GP}(A)}, \ \underline {\mathcal
{GP}(C)})$ of period $1$.
\end{thm}
\noindent{\bf Proof.} \ First, by d\'evissage each of $\underline
{\mathcal {GP}(A)}$, $\underline {\mathcal {GP}(B)}$ and $\underline
{\mathcal {GP}(C)}$ has a Serre functor. In fact, since $A$ is
Gorenstein, $\mathcal {GP}(A)$ is a resolving contravariantly finite
subcategory of $A$-{\rm mod} ([EJ, Thm. 11.5.1]; also [AR, Prop.
5.1]), and hence $\mathcal{GP}(A)$ is a resolving functorially
finite subcategory of $A$-{\rm mod} ([KS, Corol. 0.3]). Then by [AS,
Thm. 2.4] $\mathcal {GP}(A)$ has relative Auslander-Reiten
sequences. While $\mathcal {GP}(A)$ is a {\rm Frobenius} category,
by a direct argument we see that $\underline {\mathcal {GP}(A)}$ has
Auslander-Reiten triangles, and hence by [RV, Thm. I 2.4]
$\underline {\mathcal {GP}(A)}$ has a Serre functor.

Second, there is a left recollement \begin{center}
\begin{picture}(110,33)
\put(0,20){\makebox(-22,1){$\underline {\mathcal {GP}(B)}$}}
\put(40,26){\vector(-1,0){30}}
\put(10,17){\vector(1,0){30}}

\put(50,20){\makebox(25,0.8){$\underline {\mathcal {GP}(A)}$}}
\put(115,26){\vector(-1,0){30}}
\put(85,17){\vector(1,0){30}}

\put(125,20){\makebox(27,0.5){$\underline {\mathcal {GP}(C)}$}}
\put(25,30){\makebox(3,1){\scriptsize$i^\ast$}}
\put(25,21){\makebox(3,1){\scriptsize$i_\ast$}}
\put(100,30){\makebox(3,1){\scriptsize$j_!$}}
\put(100,21){\makebox(3,1){\scriptsize$j^\ast$}}
\end{picture}
\end{center}
\vskip-10pt\noindent In fact,  $_CM_B$ is compatible ([Z, Thm.
2.2{\rm(iv)}]), and hence by {\rm [Z, Thm. 1.4]}, an $A$-module \
$(X_B, Y_C)_\phi$ is in $\mathcal {GP}(A)$ if and only if $\phi:
Y\otimes_CM\rightarrow X$ is injective, ${\rm Coker}\phi\in \mathcal
{GP}(B)$, and $Y\in \mathcal {GP}(C)$. So by {\rm [Z, Thm. 3.3]} we
get the left recollement above, where $i^*$ sends $(X, Y)_\phi$ to
${\rm \Cok}\phi,$ \ \ $i_*$ sends $X$ to $(X, 0),$ \ \ $j_!$ sends
$Y$ to $(Y\otimes_CM, Y)_{\rm Id},$ and $j^*$ sends $(X, Y)_{\phi}$
to $Y$.

Now the assertion follows from Proposition
\ref{upperrecollwithserre}$(1){\rm(i)}$. $\s$

\subsection{} Recollement $(1.1)$ is {\it splitting}, if $i^! \cong
i^*$ and $j_* \cong j_!$. A splitting
recollement clearly induces a ladder of period $1$.
The product $\mathcal C'\times \mathcal C''$ of triangulated
categories $(\mathcal C', \mathcal E', T')$ and $(\mathcal C'',
\mathcal E'', T'')$ is again triangulated, where the shift $T'\times
T''$ is given by $(T'\times T'')(C', C''): = (T'C', T''C'')$, and
$\mathcal E'\times \mathcal E''$ is the collection of triangles of
$\mathcal C'\times \mathcal C''$ of the form $(X',
X'')\stackrel{(u', u'')}\longrightarrow (Y', Y'')\stackrel{(v',
v'')} \longrightarrow (Z', Z'')\stackrel{(w', w'')}\longrightarrow
(T'X', T''X''),$ where $X'\stackrel{u'}\longrightarrow
Y'\stackrel{v'} \longrightarrow Z'\stackrel{w'}\longrightarrow T'X'$
belongs to $\mathcal E'$, and $X''\stackrel{u''}\longrightarrow
Y''\stackrel{v''} \longrightarrow Z''\stackrel{w''}\longrightarrow
T''X''$ belongs to $\mathcal E''$. Then $(\mathcal C', \mathcal
C'\times \mathcal C'', \mathcal C'', p_1, \sigma_1, p_1, \sigma_2,
p_2, \sigma_2)$ is a splitting recollement, where $p_1$ and $ p_2$
are the projections, and $\sigma_1$ and $ \sigma_2$ are the
embeddings. As we see below, this gives all the splitting
recollements, up to equivalences.

\begin{prop} \label{recwithi^* = i^!} \ Let $(\mathcal C',
\mathcal C, \mathcal C'', i^\ast, i_\ast, i^!, j_!, j^\ast, j_\ast)$
be a recollement of triangulated categories. Then the following are
equivalent$:$

${\rm(i)}$ \ \   it is splitting$;$

${\rm(ii)}$ \ \   $i^! \cong i^*;$

${\rm(iii)}$ \  $j_*\cong j_!;$

${\rm(iv)}$ \  There is an equivalence $({\rm Id}_{\mathcal C'}, \
F, \ {\rm Id}_{\mathcal C''}): (\mathcal C', \mathcal C, \mathcal
C'') \longrightarrow (\mathcal C', \mathcal C'\times \mathcal C'',
\mathcal C'')$ of recollements.
\end{prop}
{\it A stable $t$-structure} ([M]) on triangulated category $\mathcal
C$  is a pair $(\mathcal U, \mathcal
V)$ of triangulated subcategories such that it is a $t$-structure ([BBD]), i.e., $\Hom(\mathcal U,
\mathcal V) =0,$ and for $X\in\mathcal C$ there is a distinguished
triangle $U\rightarrow X \rightarrow V\rightarrow U[1]$ with
$U\in\mathcal U$ and $V\in\mathcal V$. We call this triangle {\it
the $t$-decomposition} of $X$, and $U$ and $V$ {\it the $t$-part}
and {\it the $t$-free part} of $X$, respectively.
\begin{lem}\label{ffadj} \ $(1)$ \ {\rm ([M], [IKM])} \ ${\rm(i)}$ \  Given a diagram of triangle functors
\begin{picture}(70,20)
\put(20,3){\makebox(-22,1) {$\mathcal C'$}}
\put(50,8){\vector(-1,0){30}} \put(50,3){\makebox(25,0.8){$\mathcal
C$}} \put(20,0){\vector(1,0){30}}
\put(35,12){\makebox(3,1){\scriptsize$i^\ast$}}
\put(35,3){\makebox(3,1){\scriptsize$i_\ast$}}
\end{picture}
such that  $(i^{*}, i_{*})$ is an adjoint pair and $i_*$ is fully
faithful, then $({\rm Ker} i^*, \ {\rm Im}i_{*})$ is a stable
$t$-structure on $\mathcal C$, and $Y\rightarrow X
\xrightarrow{\eta_{_X}} i_*i^*X\rightarrow Y[1]$ is the
$t$-decomposition of $X$, where $\eta: {\rm Id}_{\mathcal
C}\rightarrow i_*i^*$ is the unit.

${\rm(ii)}$ \  Given a diagram of triangle functors
\begin{picture}(70,20)
\put(20,4){\makebox(-22,1) {$\mathcal C'$}}
\put(50,-1){\vector(-1,0){30}} \put(50,4){\makebox(25,0.8){$\mathcal
C$}} \put(20,8){\vector(1,0){30}}
\put(35,3){\makebox(3,1){\scriptsize$i^!$}}
\put(35,12){\makebox(3,1){\scriptsize$i_\ast$}}
\end{picture}
such that $(i_{*}, i^!)$ is an adjoint pair and $i_*$ is fully
faithful, then $({\rm Im}i_{*}, \ {\rm Ker} i^!)$ is a stable
$t$-structure on $\mathcal C$, and $i_*i^!X \xrightarrow
{\epsilon_{_X}} X\rightarrow Z\rightarrow (i_*i^!X)[1]$ is the
$t$-decomposition of $X$, where $\epsilon: i_*i^!\rightarrow {\rm
Id}_{\mathcal C}$ is the counit.

\vskip5pt

$(2)$ \ Let $(\mathcal Y,
\mathcal Z)$ be a stable $t$-structure on \ $\mathcal C$ with
$\Hom(\mathcal Z, \mathcal Y) = 0$. Then $F: \mathcal C\rightarrow
\mathcal Y \times \mathcal Z$  given by $FX = (Y, Z)$ is a
triangle-equivalence, where $Y\stackrel u\rightarrow  X \rightarrow
Z \rightarrow Y[1]$ is the $t$-decomposition.
\end{lem}
\noindent{\bf Proof.} \ $(2)$ \ By assumption $\Hom_\mathcal
C(Z[-1], Y) = 0$. By the exact sequence $\Hom_\mathcal C(X, Y)
\stackrel{\Hom(u, Y)}\longrightarrow \Hom_\mathcal C(Y, Y) $
$\rightarrow\Hom_\mathcal C(Z[-1], Y) = 0$ we see that $u$ is a
splitting monomorphism. Thus $X\cong Y\oplus Z$ ([H1, p.7]). It is
straightforward that $F: \mathcal C\rightarrow \mathcal Y \times
\mathcal Z$ given by $FX = (Y, Z)$ is a triangle-equivalence. $\s$

\vskip5pt

\noindent{\bf Proof of Proposition \ref{recwithi^* = i^!}.}
${\rm(i)}\Longrightarrow {\rm(ii)}$ and ${\rm(iv)}\Longrightarrow
{\rm(i)}$ are obvious.

${\rm(ii)}\Longrightarrow {\rm(iii)}:$ \ \ Suppose $i^! \cong i^*.$
For $X\in\mathcal C$ and $Y''\in\mathcal C''$ applying
$\Hom_\mathcal C(-, j_!Y'')$ to the recollement triangle $j_!j^*X
\rightarrow X \rightarrow i_*i^*X \rightarrow
(j_!j^*X)[1]$ we get the exact sequence
$$\Hom(i_*i^*X, j_!Y'')\longrightarrow \Hom(X, j_!Y'')\longrightarrow \Hom(j_!j^*X, j_!Y'')\longrightarrow
\Hom((i_*i^*X)[-1], j_!Y'').$$ By $\Hom(i_*i^*X,
j_!Y'')\cong\Hom(i^*X, i^!j_!Y'') \cong \Hom(i^*X, i^*j_!Y'') = 0$ and
$\Hom((i_*i^*X)[-1], j_!Y'') = 0,$ we have
$\Hom_\mathcal C(X, j_!Y'')\cong \Hom_{\mathcal C}(j_!j^*X,
j_!Y'')\cong \Hom_{\mathcal C''}(j^*X, Y'')$, i.e., $(j^*,
j_!)$ is an adjoint pair. While $(j^*, j_*)$ is also an adjoint
pair, so  $j_*\cong j_!.$

${\rm(iii)}\Longrightarrow {\rm(ii)}$ can be similarly proved.

${\rm(i)}\Longrightarrow {\rm(iv)}:$  \ \ Assume that $i^! \cong
i^*$ and $j_*\cong j_!.$ Since $(i_*, i^!)$ is an adjoint pair, so
is $(i_*, i^*)$, and hence by Lemma \ref{ffadj}$(1){\rm(ii)}$  $({\rm Im}i_*,
{\rm Ker}i^*)$ is a stable $t$-structure. Since $\Hom({\rm Ker}i^*, \ {\rm
Im}i_* ) = 0$ and the recollement
triangle $i_*i^!X \rightarrow X \rightarrow j_*j^*X
\rightarrow (i_*i^!X)[1]$ is the $t$-decomposition (since
$j_*j^*X\in {\rm Im}j_* = {\rm Ker}i^! = {\rm Ker}i^*$ by the
assumption), by Lemma \ref{ffadj}$(2)$ $\widetilde{F}:
\mathcal C\rightarrow {\rm Im}i_* \times {\rm Ker}i^*$ given by
$\widetilde{F}X = (i_*i^!X, j_*j^*X)$ is a triangle-equivalence.
Since ${\rm Im}i_*\cong \mathcal C'$ and \ ${\rm Ker}i^* = {\rm
Im}j_! \cong \mathcal C''$,  we get a
triangle-equivalence $F: \mathcal C\rightarrow \mathcal C' \times
\mathcal C''$ with $FX = (i^!X, j^*X)$. Now it is straightforward
that $({\rm Id}_{\mathcal C'},  F, {\rm Id}_{\mathcal C''}):
(\mathcal C', \mathcal C, \mathcal C'') \rightarrow (\mathcal C',
\mathcal C'\times \mathcal C'', \mathcal C'')$ is an equivalence of
recollements. We omit the details. $\s$
\begin{rem} \ ${\rm(i)}$ \ A Hom-finite $k$-triangulated category
$(\mathcal C, [1])$ is a $d$-{\it Calabi-Yau category} {\rm([Ke2])}, if
there is a nonnegative integer $d$, such that the $d$-th shift $[d]$
is a right Serre functor of $\mathcal C$.

By Lemma \ref{withserre} any left {\rm(}right{\rm)} recollement of Calabi-Yau
category $\mathcal C$ sits in  a splitting recollement. Thus any
recollement of Calabi-Yau category is splitting.

\vskip5pt

${\rm(ii)}$ \ If $(\mathcal C', \mathcal C, \mathcal C'')$ is a recollement with $\mathcal C$
{\rm Calabi-Yau}, then obviously so are $\mathcal C'$ and $\mathcal C''$. However, the converse is not true$:$
otherwise, $(\mathcal C', \mathcal C, \mathcal C'')$ is splitting by
${\rm(i)};$ but there are a lot of
examples of non-splitting recollements $(\mathcal C', \mathcal C,
\mathcal C''),$ where $\mathcal C'$ and $\mathcal C''$ are {\rm
Calabi-Yau}. For example, let $A=\left(\begin{smallmatrix}
k&0\\
k&k
\end{smallmatrix}\right)$ with $k$ a field.
Then one has a recollement $(D^b(k\mbox{-}{\rm mod}),
D^b(A\mbox{-}{\rm mod}), D^b(k\mbox{-}{\rm mod}))$ {\rm(}{\rm[PS,
Exam. 2.10]}{\rm)}. Note that $D^b(k\mbox{-}{\rm mod})$ is
$0$-{\rm Calabi-Yau} and  that $(D^b(k\mbox{-}{\rm mod}),
D^b(A\mbox{-}{\rm mod}), D^b(k\mbox{-}{\rm mod}))$ is not splitting
{\rm(}otherwise, $D^b(A\mbox{-}{\rm mod})$ is the product of two
{\rm Calabi-Yau} categories, and hence again {\rm
Calabi-Yau}{\rm}$;$ but $D^b(A\mbox{-}{\rm mod})$ is not {\rm
Calabi-Yau)}.
\end{rem}

\centerline {\bf Appendix: Proofs of lemmas in Section 1}

\vskip5pt

We include proofs of lemmas in Section 1 only for convenience
(although they are well-known, it seems that explicit proofs are not
available in the literature).

\vskip5pt

{\bf Proof of Lemma \ref{criterionupperrec}.} \  Since a right
recollement of $\mathcal C$ relative to $\mathcal C'$ and $\mathcal
C''$ is a left recollement of $\mathcal C$ relative to $\mathcal
C''$ and $\mathcal C'$,  \ it follows that $(2)$ can be deduced from
$(1)$. We include a proof of $({\rm ii})\Longrightarrow ({\rm i})$
of $(1)$.

Since $(i^*, i_*)$ is an adjoint pair and $i_*$ is fully faithful, by Lemma \ref{ffadj}$(1)$${\rm(i)}$ $Y\rightarrow X \xrightarrow{\eta_{_X}} i_*i^*X\rightarrow Y[1]$
is the $t$-decomposition of $X$ respect to the $t$-structure $({\rm Ker}i^*, \ {\rm Im}i_*).$
Similarly, by Lemma \ref{ffadj}$(1)$${\rm(ii)}$ $j_!j^*X \xrightarrow {\epsilon_{_X}} X\rightarrow Z\rightarrow (j_!j^*X)[1]$
is the $t$-decomposition of $X$ respect to the $t$-structure $({\rm Im}j_!, \ {\rm Ker}j^*).$
Since both $({\rm Ker}i^*, \ {\rm Im}i_*)$ and $({\rm Im}j_!, \ {\rm Ker}j^*)$ are $t$-structures and
${\rm Im}i_\ast = {\rm Ker}j^*$, it follows that ${\rm Ker}i^* = {\rm Im}j_!$, and the two $t$-decompositions above are isomorphic. From this one easily deduces that
$j_!j^*X
\stackrel{\epsilon_{_X}}\longrightarrow X
\stackrel{\eta_{_X}}\longrightarrow i_*i^*X\rightarrow
(j_!j^*X)[1]$ is a distinguished triangle. $\s$

\vskip5pt

\noindent{\bf Lemma A.1.} \label{adjquo} {\rm (see e.g. [BBD], [M],
[N3], [IKM])} \ {\it Let $(\mathcal U, \mathcal V)$ be a stable
$t$-structure  on $\mathcal C$. Then

${\rm(i)}$ \ there is a triangle-equivalence $V_\mathcal
V\circ\sigma_\mathcal U: \mathcal U\rightarrow\mathcal C/\mathcal
V$, where $\sigma_\mathcal U: \mathcal U\hookrightarrow\mathcal C$
is the embedding, and $V_\mathcal V: \mathcal C\rightarrow\mathcal
C/\mathcal V$ is the Verdier functor. A quasi-inverse of $V_\mathcal
V\circ\sigma_\mathcal U$ sends object $X\in \mathcal C/\mathcal V$
to its $t$-part.

${\rm(ii)}$ \ there is a triangle-equivalence $V_\mathcal
U\circ\sigma_\mathcal V: \mathcal V\rightarrow\mathcal C/\mathcal
U$, where $\sigma_\mathcal V: \mathcal V\hookrightarrow\mathcal C$
is the embedding, and $V_\mathcal U: \mathcal C\rightarrow\mathcal
C/\mathcal U$ is the Verdier functor. A quasi-inverse of $V_\mathcal
U\circ\sigma_\mathcal V$ sends object $X\in \mathcal C/\mathcal U$
to its $t$-free part.}

\vskip5pt

\noindent{\bf Lemma A.2.} \label{rladjandfaithful} {\rm([AHKLY, Lemma 2.2])} \ {\it Let \ $\mathcal C'\stackrel F\longrightarrow\mathcal C\stackrel G\longrightarrow\mathcal C''$ \ be a sequence of triangle functors, such that
$F$ is fully faithful, ${\rm Im} F = {\rm Ker}G$, and $G$ induces a triangle-equivalence $\mathcal C/{\rm Ker}G\cong\mathcal C''.$ Then
$F$ has a right {\rm(}resp. left{\rm)} adjoint $F'$ if and only if so does $G$.

In this case, the right {\rm(}resp. left{\rm)} adjoint $G'$ of $G$ is also fully faithful,
${\rm Im} G' = {\rm Ker}F'$, and $F'$ induces a triangle-equivalence $\mathcal C/{\rm Ker}F'\cong\mathcal C'.$}

\vskip5pt

\noindent{\bf Proof.} \ Using the opposite category, we only need to prove the right version.

By the universal property, $G$ is the composition of the Verdier
functor $\mathcal C\longrightarrow\mathcal C/{\rm Ker}G$ with the
equivalence $\mathcal C/{\rm Ker}G\cong\mathcal C''$. Thus, for
simplicity, without loss of the generality we may assume that
$\mathcal C/{\rm Ker}G=\mathcal C''$ and $G$ is just the Verdier
functor $\mathcal C\rightarrow\mathcal C/{\rm Ker}G.$

$\Longleftarrow:$ \ Assume that $G$ has a right adjoint pair $G'$,
i.e., a Bousefield localization functor exists for the pair ${\rm
Ker}G\subseteq \mathcal C$. Thus for $X\in\mathcal C$, by A. Neeman
[N3, Prop. 9.1.8] there is a distinguished triangle $Z\rightarrow
X\stackrel{\eta_X}\longrightarrow G'GX\rightarrow Z[1]$ with $Z\in
{\rm Ker}G ={\rm Im}F$, where $\eta: {\rm Id}_\mathcal C \rightarrow
G'G$ is the unit. Thus \ $({\rm Im}F, \ {\rm Im}G')$ is a
$t$-structure on $\mathcal C$, which induces an adjoint pair
$(\sigma, \widetilde{F'})$, where $\sigma: {\rm Im}F\rightarrow
\mathcal C$ is the embedding, and $\widetilde{F'}: \mathcal
C\rightarrow {\rm Im}F$ sends $X$ to its $t$-part $Z$. Since $Z\in
{\rm Im}F$ and $F$ is fully faithful, there is a unique object (up
to isomorphism) $Z'\in \mathcal C'$ such that $Z \cong FZ'$. Define
$F': \ \mathcal C\rightarrow\mathcal C'$ to be the functor given by
$X\mapsto Z'$. Since $(\sigma, \widetilde{F'})$ is an adjoint pair
and $F$ is fully faithful, it is easy to see that $(F, F')$ is an
adjoint pair. By construction we have ${\rm Im} G' = {\rm Ker}F'.$
Since $({\rm Im}F, \ {\rm Im}G')$ is a $t$-structure, it follows
from Lemma A.1${\rm(i)}$ that $X\mapsto Z$ gives an
triangle-equivalence $\mathcal C/{\rm Im} G'\rightarrow {\rm Im}F;$
together with ${\rm Im}F\cong \mathcal C'$ we see that $F'$ induces
a triangle-equivalence $\mathcal C/{\rm Ker}F'\cong\mathcal C'.$
Since $G(Z) = 0$, $G(\eta_X)$ is an isomorphism, and hence by
$\epsilon_{GX}\circ G(\eta_X) = {\rm Id}_{\mathcal C''}$ (where
$\epsilon: GG'\rightarrow {\rm Id}_{\mathcal C''}$ is the counit) we
see that $\epsilon_{GX}$ is an isomorphism for each $X\in\mathcal
C$. Since by assumption $G$ is dense, $\epsilon: GG'\rightarrow {\rm
Id}_{\mathcal C''}$ is a natural isomorphism of functors, and thus
$G'$ is fully faithful.

$\Longrightarrow:$ \ Assume that $F$ has a right adjoint pair $F'$.
Then by Lemma \ref{ffadj}$(1)$${\rm(ii)}$ $({\rm Im}F, \ {\rm
Ker}F')$ is a $t$-structure on $\mathcal C$, with $t$-decomposition
$FF'X\stackrel{\omega_X}\longrightarrow X\rightarrow Y\rightarrow
(FX')[1]$  of $X\in \mathcal C$, where $\omega: FF'\rightarrow {\rm
Id}_\mathcal C$ is the counit. This $t$-structure induces an adjoint
pair $(\widetilde{G}, \sigma)$, where $\widetilde{G}: \mathcal
C\rightarrow {\rm Ker}F'$ sends $X$ to its $t$-free part $Y$, and
$\sigma: {\rm Ker}F'\rightarrow \mathcal C$ is the embedding. By
Lemma A.1${\rm(ii)}$ the functor $\widetilde{G'}: \ \mathcal C/{\rm
Im}F\rightarrow {\rm Ker}F'$, which sends each object $X$ to its
$t$-free part $Y$, is a triangle-equivalence. Thus $G =
\widetilde{G'}^{-1}\widetilde{G}$. Put $G': =\sigma\widetilde{G'}:
\mathcal C/{\rm Im}F\rightarrow \mathcal C$, i.e., \ $G': \mathcal
C''\rightarrow\mathcal C$. By construction $G'$ is fully faithful
and ${\rm Im} G' = {\rm Ker}F'.$ By Lemma A.1${\rm(i)}$ $\mathcal
C/{\rm Ker} F'\rightarrow {\rm Im}F$ given by $X\mapsto FF'X$ is an
triangle-equivalence; together with ${\rm Im}F\cong \mathcal C'$ we
see that $F'$ induces $\mathcal C/{\rm Ker}F'\cong\mathcal C'.$ For
$X\in\mathcal C$ and $C''\in\mathcal C''$,  since $(\widetilde{G},
\sigma)$ is an adjoint pair, we have
$$\Hom(GX, C'') = \Hom(\widetilde{G'}^{-1}\widetilde{G}X, C'') \cong
\Hom_{{\rm Ker}F'}(\widetilde{G}X, \widetilde{G'}C'')\cong \Hom_{\mathcal C}(X, \sigma\widetilde{G'}C'')=\Hom(X, G'C''),$$
i.e., $(G, G')$ is an adjoint pair. $\s$

\vskip5pt

{\bf Proof of Lemma \ref{laddercriterion}.} \ It suffices to prove the ``if" part. We denote the recollement ${\rm (1.1)}$ by $(\mathcal C', \mathcal C, \mathcal
C'',$ $j_{-1}, i_0, j_1, i_{-1}, j_0, i_1)$ (this labeling coincides with $(1.2)$), and assume that
there is an infinite adjoint sequence \ ${\rm(}\cdots, i_{-2}, \ j_{-1}, i_0, j_1, i_2, \cdots{\rm)}$.
Since $i_1$ is fully faithful and $j_1$ has a right adjoint pair $i_2$, by applying Lemma A.2 to the sequence
$\mathcal C''\stackrel {i_1}\longrightarrow\mathcal C\stackrel {j_1}\longrightarrow\mathcal C'$ we get an adjoint pair $(i_1, j_2)$,  such that
the right adjoint of $j_1$ is fully faithful (i.e., $i_2$ is fully faithful),
${\rm Im} i_2 = {\rm Ker}j_2$, and that $j_2$ induces a triangle-equivalence $\mathcal C/{\rm Ker}j_2\cong\mathcal C''.$ Applying Lemma A.2 to the sequence
$\mathcal C'\stackrel {i_2}\longrightarrow\mathcal C\stackrel {j_2}\longrightarrow\mathcal C''$, and continuing this process we then
get a ladder going downwards infinitely, by Lemma \ref{criterionupperrec}.

Going upwards, and by the same argument we get a ladder going upwards infinitely. Putting together we get an unbounded ladder containing recollement $(\mathcal C', \mathcal C, \mathcal
C'', j_{-1}, i_0, j_1, i_{-1}, j_0, i_1)$. $\s$

\vskip10pt

{\bf Proof of Lemma \ref{periodicladder}.} \ $(1)$ \ We only prove
${\rm (ii)} \Longrightarrow {\rm (i)}.$ \ Any recollement $(\mathcal
C', \mathcal C, \mathcal C'', i^*, i_*, i^!, j_!, j^*, j_*)$ induces
an equivalence $(\widetilde{i_*}, {\rm Id}_\mathcal
C, \widetilde{j_*}): (\mathcal C', \mathcal C, \mathcal C'', i^*,
i_*, i^!, j_!, j^*, j_*)\rightarrow ({\rm Im}i_*, \ \mathcal C,
\  {\rm Im}j_*, \ \widetilde{i_*}i^\ast, \ \sigma_1, \
\widetilde{i_*}i^!, \ \widetilde{j_!}, \ \widetilde{j_*}j^\ast, \
\sigma_2)$ of recollements, where $\widetilde{i_*}: \mathcal C'\rightarrow {\rm
Im}i_*$ and $\widetilde{j_*}: \mathcal C''\rightarrow {\rm
Im}j_*$ are the equivalences induced by $i_*$ and $j_*$, respectively, $\sigma_1$ and $\sigma_2$ are
embeddings, and $\widetilde{j_!}: {\rm Im}j_*\rightarrow
\mathcal C$ is given by $j_*C''\mapsto j_!C'', \ \forall \
C''\in\mathcal C''.$ By restriction we get $\widetilde{F'}: {\rm
Im}i_*\stackrel\sim\longrightarrow {\rm Im}i^\mathcal D_*$ and
$\widetilde{F''}: {\rm Im}j_*\stackrel\sim\longrightarrow {\rm
Im}j^\mathcal D_*.$  Thus, it suffices to prove that there is an equivalence
\vskip5pt
\begin{center}
\begin{picture}(100,100)
\put(-2,80){\makebox(3,1){${\rm Im}i_*$}}
\put(40,90){\vector(-1,0){30}}
\put(10,80){\vector(1,0){30}}
\put(40,70){\vector(-1,0){30}}

\put(50,80){\makebox(3,1){${\mathcal C}$}}

\put(90,90){\vector(-1,0){30}}
\put(60,80){\vector(1,0){30}}
\put(90,70){\vector(-1,0){30}}
\put(102,80){\makebox(3,1){${\rm
Im}j_*$}}

\put(25,95){\makebox(3,1){\scriptsize$\widetilde{i_*}i^\ast$}}
\put(25,85){\makebox(3,1){\scriptsize$\sigma_1$}}
\put(25,75){\makebox(3,1){\scriptsize$\widetilde{i_*}i^!$}}

\put(75,95){\makebox(3,1){\scriptsize$\widetilde{j_!}$}}
\put(75,85){\makebox(3,1){\scriptsize$\widetilde{j_*}j^\ast$}}
\put(75,75){\makebox(3,1){\scriptsize$\sigma_2$}}

\put(-12,50){$\widetilde{F'}$}

\put(40,50){$F$}

\put(105,50){$\widetilde{F''}$}

\put(0,70){\vector(0,-1){40}}

\put(50,70){\vector(0,-1){40}} \put(100,70){\vector(0,-1){40}}

\put(-3,25){\makebox(3,1){${\rm Im}i_*^\mathcal D$}}
\put(40,35){\vector(-1,0){30}} \put(10,23){\vector(1,0){30}}
\put(40,13){\vector(-1,0){30}} \put(50,25){\makebox(3,1){${\mathcal
D}$}} \put(90,35){\vector(-1,0){30}} \put(60,23){\vector(1,0){30}}
\put(90,13){\vector(-1,0){30}} \put(104,25){\makebox(3,1){${\rm
Im}j_*^\mathcal D$}}

\put(25,40){\makebox(3,1){\scriptsize$\widetilde{i_*^\mathcal
D}i_D^\ast$}} \put(25,28){\makebox(3,1){\scriptsize$\sigma_1^D$}}
\put(25,17){\makebox(3,1){\scriptsize$\widetilde{i_*^\mathcal D}
i_\mathcal D^!$}}
\put(75,40){\makebox(3,1){\scriptsize$\widetilde{j^\mathcal D_!}$}}
\put(75,28){\makebox(3,1){\scriptsize$\widetilde{j_*^\mathcal D}
j_\mathcal D^\ast$}}
\put(75,18){\makebox(3,1){\scriptsize$\sigma_2^\mathcal D$}}
\end{picture}
\end{center}
\vskip-10pt\noindent  i.e.,
for
$C\in\mathcal C$ and $j_*C''\in {\rm Im}j_*$ with $C''\in\mathcal
C''$, there are natural isomorphisms:
$F{i_*}i^\ast C \cong {i_*^\mathcal D} i^\ast_\mathcal D FC,
\ \ F{i_*}i^!C\cong {i_*^\mathcal D} i^!_\mathcal D FC, \ \
Fj_!C''\cong \widetilde{j^\mathcal D_!} Fj_*C'', \ \ F {j_\ast}
j^*C\cong {j^\mathcal D_\ast} j^*_\mathcal D FC.$ By the
recollement triangle $i_*i^!C \rightarrow C \rightarrow
j_*j^*C\rightarrow (i_*i^!C)[1]$ we get distinguished triangles
$$Fi_*i^!C \rightarrow FC \rightarrow Fj_*j^*C\rightarrow (Fi_*i^!C)[1], \ \
\mbox{and} \ \ i_*^\mathcal Di^!_\mathcal DFC \rightarrow FC
\rightarrow j_*^\mathcal Dj^*_\mathcal DFC\rightarrow (i^\mathcal
D_*i_\mathcal D^!FC)[1].$$ By the assumption, they are both the $t$-decompositions of $FC$ respect to the $t$-structure
$({\rm Im}i_*^\mathcal D, {\rm Im}j^\mathcal D_*)$, hence
$F{i_*}i^!C\cong {i_*^\mathcal D }i^!_\mathcal D FC$ and
$F{j_*}j^*C\cong j_*^\mathcal D j^*_\mathcal D FC$. Similarly, by $j_!j^*C \rightarrow C
\rightarrow i_*i^*C\rightarrow (j_!j^*C)[1]$  we get $Fj_!j^*C \cong j_!^\mathcal Dj^*_\mathcal
DFC$ and $Fi_*i^*C\cong i_*^\mathcal Di^*_\mathcal DFC$.
It remains to prove $Fj_!C''\cong \widetilde{j^\mathcal D_!}
Fj_*C''$.  By $C''\cong j^*j_*C''$ the functor $\widetilde{j_!}$
reads as $\widetilde{j_!} j_*C''= j_!j^*j_*C''.$ Since $Fj_*C''\in
{\rm Im}j_*^\mathcal D$, we have $\widetilde{j^\mathcal
D_!}Fj_*C''\cong j_!^\mathcal Dj^*_\mathcal D Fj_*C''.$ It follows
that $Fj_!C'' \cong Fj_!j^*j_*C''\cong j_!^\mathcal Dj^*_\mathcal
DFj_*C''\cong \widetilde{j^\mathcal D_!}Fj_*C''.$

\vskip5pt

${\rm (2)}$ \ We claim that the $t$-th recollement is equivalent to
the $0$-th one. In fact, by assumption there is equivalence $(F', F,
F''): (\mathcal C', \mathcal C, \mathcal C'', j_{2t-1}, i_{2t},
i_{2t-1}, j_{2t})\rightarrow (\mathcal C', \mathcal C, \mathcal C'',
j_{-1}, i_0, i_{-1}, j_0)$ of left recollements. It remains to prove
that there  are natural isomorphisms $F'j_{2t+1}\cong j_1F$ and
$Fi_{2t+1}\cong i_1F''$. Since $(i_{2t}, \ j_{2t+1})$ and $(j_{2t},
\ i_{2t+1})$ are adjoint pairs, it suffices to prove $(i_{2t},
F'^{-1}j_1F)$ and $(j_{2t}, F^{-1}i_1F'')$ are also adjoint pairs.
Indeed, the first adjoint pair can be seen from (and the second one
is similarly proved)
$$\Hom(X', F'^{-1}j_1FY)\cong \Hom(F'X',
j_1FY) \cong \Hom(i_0F'X', FY)\cong
\Hom(Fi_{2t}X', FY) \cong \Hom(i_{2t}X',
Y).$$

Going downwards (resp. upwards) step by step, by the similar argument we
see the assertion.

\vskip5pt

${\rm (3)}$ follows from ${\rm (1)}$ and $(2)$. $\s$

\vskip10pt

\noindent {\small \it Pu Zhang, College of Mathematics, \  Shanghai
Jiao Tong University, \ Shanghai 200240, China,  \
pzhang$\symbol{64}$sjtu.edu.cn}

\vskip5pt

\noindent {\small\it Yuehui Zhang, College of Mathematics, \
Shanghai Jiao Tong University,  \ Shanghai 200240, China, \
zyh$\symbol{64}$sjtu.edu.cn}

\vskip5pt

\noindent {\small\it Guodong Zhou, Department of Mathematics, \
Shanghai Key laboratory of PMMP, \ East China Normal University,  \
Shanghai 200241, China, \ gdzhou$\symbol{64}$math.ecnu.edu.cn}

\vskip5pt

\noindent {\small\it Lin Zhu, College of Mathematics, \  Shanghai
Jiao Tong University,  \ Shanghai 200240, China, \
zhulin2323$\symbol{64}$163.com}

\end{document}